\documentclass{amsart}

\usepackage{amssymb}
\usepackage{hyperref}
\usepackage{color}
\usepackage{fancyhdr}
\usepackage{framed}
\usepackage{enumerate}

\newtheorem{Theorem}{Theorem}
\newtheorem{Definition}[Theorem]{Definition}
\newtheorem{Proposition}[Theorem]{Proposition}

\newtheorem{Lemma}[Theorem]{Lemma}
\newtheorem{Corollary}[Theorem]{Corollary}
\theoremstyle{remark}

\newtheorem{Example}[Theorem]{Example}

\def\il{\int_}
\def\hvi{\varphi}
\def\eps{\varepsilon}

\def\ovr{\overline}

\def\al{\alpha}

\def\th{\theta}

\def\Dl{\Delta}
\def\dl{\delta}

\def\bd{\partial}
\def\lm{\lambda}

\def\sm{\setminus}
\def\sbs{\subset}

\def\nea{\nearrow}
\def\sea{\searrow}

\def\be{\begin{enumerate}}
\def\ee{\end{enumerate}}
\def\bT{\begin{Theorem}}
\def\eT{\end{Theorem}}
\def\bP{\begin{Proposition}}
\def\eP{\end{Proposition}}
\def\bD{\begin{Definition}}
\def\eD{\end{Definition}}
\def\bE{\begin{Example}}
\def\eE{\end{Example}}
\def\bL{\begin{Lemma}}
\def\eL{\end{Lemma}}
\def\bC{\begin{Corollary}}
\def\eC{\end{Corollary}}

\def\oD{\ovr{\mathbb D}}
\def\aD{\mathbb D}

\def\aT{\mathbb T}
\def\aC{\mathbb C}

\def\E{{\mathcal E}}

\begin{document}
\title{Weighted Hardy spaces on the unit disk}
\author{Khim R. Shrestha}
\address{Department of Mathematics,  Syracuse University, \newline
215 Carnegie Hall, Syracuse, NY 13244} \email{krshrest@syr.edu}
\maketitle
\begin{abstract}
In this paper we mainly discuss three things. First, there is no canonical norm on the space $H^p_u(\aD)$. Second, we improve the weak-$*$ convergence of the measures $\mu_{u,r}$. Third, the dilations $f_t$ of the function $f\in H^p_u(\aD)$ converge to $f$ in $H^p_u$-norm and hence the polynomials are dense in $H^p_u(\aD)$.  
\end{abstract}
\section{Introduction}
\par Suppose that $u$ is a continuous function on the closure of a domain $D\sbs\aC^n$ that is plurisubharmonic on $D$ and equal to 0 on $\bd D$. In \cite{D1} for each $r<0$ Demailly introduced a measure $\mu_{u,r}$ supported by the set $S_{u,r}=\{u=r\}$  (see Section \ref{basic} for precise definition). Using these measures Poletsky and Stessin introduced in \cite{PS} the weighted Hardy spaces on $D$ as the space $H^p_u(D)$ of all holomorphic functions $f$ on $D$ such that
\[\|f\|_{H^p_u}^p = \varlimsup_{r\to 0^-}\il {S_{u,r}} |f|^p\, d\mu_{u,r} < \infty.\]
In the same manner the weighted Hardy spaces $h^p_u(D)$ of pluriharmonic functions also can be introduced. 
\par Recently these spaces were subjected to more detailed studies in \cite{MN, Sh, Sibel} when $D$ is the unit disk $\aD$. The main goal of these papers was to establish properties of weighted Hardy spaces similar to the standard  properties of the classical Hardy spaces. In this paper we continue this program looking more thoroughly at the properties of boundary values.
\par It was proved in \cite{PS} that $H^p_u(D)\sbs H^p(D)$ for all exhausting functions $u$. Hence every function $f$ in $H^p_u(\aD)$ or in $h^p_u(\aD)$ has radial boundary values $f^*$ almost everywhere on $\aT=\bd\aD$ with respect to the Lebesgue measure. It was proved in \cite{D1} that the measures $\mu_{u,r}$ converge weak-$*$ in $C^*(\ovr D)$ to a measure $\mu_u$ supported by $\aT$ when $r\to0^-$ and it was established in \cite{MN, Sh, Sibel} that
\[\|f\|_{H^p_u}=\|f^*\|_{L^p_u}\]
where $L^p_u(\aT) = L^p(\aT, \mu_u)$.
\par In Section \ref{basic} we list all necessary definitions and known facts. In Section \ref{norm} we show that the $h^p_u$-norm of an $h^p_u$-function is equal to the $L^p_u$-norm of its boundary value function. We also observe that while different exhaustion functions may define the same space it is impossible to select a canonical exhaustion  determining a canonical norm. For example, if $\E_0$ denotes the set of all exhaustions such that the measure $(dd^cu)^n$ has the total mass equal to 1 and a compact support then it was proved in \cite{PS} that all Hardy spaces $H^p_u(D)$ coincide when $u\in\E_0$. But as we show the intersections of unit balls in these spaces is the unit ball in $H^\infty(\aD)$ when $D = \aD$. 
\par In Section \ref{wc} we show that for any function $h \in h^p_u(\aD),\, p > 1,$ the measures $\{h\mu_{u,r}\}$ converge weak-${\ast}$ to $h^*\mu_u$ in $C^*(\oD)$ when $r\to0^-$. In particular, the measures $\mu_{u,r}$ converge weak$*$ to $\mu_u$ in the dual of $h^p_u(\aD)$ when $r\to0^-$. This improves the result of Demailly mentioned above and shows that the convergence is much stronger. 
\par In Section \ref{cd} we look at convergence of the dilations $f_t(z)=f(tz)$, $0<t<1$, to the function $f\in H^p_u(\aD)$, $p > 0$. In the classical case it is known that the functions $f_t$ converge to $f$ in $H^p(\aD)$. We generalize this result to space $H^p_u(\aD)$ and this allows us to prove that polynomials are dense in $H^p_u(\aD)$. 
\par I would like to express my sincere gratitude to my advisor Prof. E. A. Poletsky for his continuous support and guidance. Without his advice this paper would never have gotten into this form. 

\section{Basic facts} \label{basic}

\par Let $\aD$ be the unit disc $\{|z|< 1\}$ in $\aC$.  A continuous subharmonic function $u:\mathbb{D}\to [-\infty, 0)$ such that $u(z)\to 0$ as $|z|\to 1$ is called an exhaustion function. Following \cite{D1} for $ r < 0 $ we set
\[B_{u,r}  =\{z\in \mathbb{D}:u(z) < r\}\text{ and }
S_{u,r} = \{z\in \mathbb{D}: u(z) = r\}.\]
As in \cite{D1} we let $u_r = \max\{u,r\}$ and define the  measure
\[\mu_{u,r} = \Dl u_r -\chi_{\aD \setminus B_r}\Dl u, \] where $\Dl$ is the Laplace operator.
Clearly $\mu_{u,r}\ge 0$ and is  supported by $S_{u,r}$.
\par Let us denote by $\E$ the set of all continuous negative subharmonic exhaustion functions $u$ on $\aD$ such that
\[\il \aD\Dl u = 1.\] In the same paper Demailly (see Theorems 1.7 and 3.1 there) proved the following result which we adapt to the case of $\aD$.
\bT [Lelong--Jensen formula]\label{jlf}
Let $\phi$ be a subharmonic function on $\mathbb{D}$. Then $\phi$ is $\mu_{u,r}$-integrable for every $r < 0$ and \[\mu_{u,r}(\phi) = \int_{B_{u,r}} \phi\, \Dl u + \int_{B_{u,r}}(r-u)\,\Dl\phi.\]
Moreover, if $u\in\E$ then the measures $\mu_{u,r}$ converge weak-$*$ in $C^*(\oD)$  to a measure $\mu_u\ge 0$ supported by $\aT$ as $r\to 0^-$.
\eT
As a consequence of this theorem he derived the following.
\bC\label{C:imu} If $\phi$ is a non-negative subharmonic function, then the function $r\mapsto \mu_{u,r}(\phi)$ is increasing on $(-\infty, 0)$.\eC
\par Using the measures $\mu_{u,r}$, in \cite{PS}, Poletsky and Stessin introduced the weighted Hardy spaces associated with an exhaustion $u\in\E$. Following \cite{PS} we define the space $H^p_u(\aD), 0 < p < \infty,$ consisting of all holomorphic functions $f(z)$ in $\mathbb{D}$ that satisfy \[\|f\|_{H^p_u}^p = \varlimsup_{r\to 0^-}\il {S_{u,r}} |f|^p\, d\mu_{u,r} < \infty.\] By Corollary \ref{C:imu} we can replace the $\varlimsup$ in the above definition with $\lim$.  By Theorem \ref{jlf} and the monotone convergence theorem it follows that,
\begin{equation}\label{nf} \|f\|^p_{H^p_u} = \il {\aD} |f|^p\,\Dl u - \il {\aD} u\,\Dl |f|^p.
\end{equation}
The classical Hardy spaces $H^p(\aD)$ correspond to the exhaustion function $u(z)=\log|z|$ (\cite[Section 4]{PS}). Hence the classical definition of the Hardy spaces is subsumed in this new definition.

\par It is proved in \cite{PS} that:\be
\item the spaces $H^p_u(\aD)$ are Banach when $p\ge 1$ (Theorem 4.1);
\item if $v,u\in\E$ and $v\le u$ on $\aD$, then $H^p_v(\aD)\sbs H^p_u(\aD)$ and if $f\in H^p_v(\aD)$ then $\|f\|_{H^p_u}^p\le \|f\|_{H^p_v}^p$.
\ee
\par Thus by Hopf's lemma the space $H^p_u(\mathbb{D})$ is contained in the classical Hardy space $H^p(\mathbb{D})$.

\par  M. A. Alan and N. G. Go\u{g}u\c{s} in \cite{MN},  S. \c{S}ahin in \cite{Sibel}  and K. R. Shrestha in \cite{Sh}  have independently produced examples that in general these new spaces $H^p_u(\aD)$ do not coincide with the classical spaces $H^p(\aD)$. They have found an exhaustion funcion $u$ for which $H^p_u(\aD) \subsetneq H^p(\aD)$.
 \par It has been established (see \cite{MN}, \cite{Sh} and \cite{Sibel}) that the boundary measure $\mu_u$ is absolutely continuous with respect to the Lebesgue measure $\lm$ on $\aT$. Here and throughout this paper $\lm$ is normalized, i.e. $\int_\aT d\lm = 1$.  Hence \[d\mu_u = \al_u d\lm\] for some $\al_u \in L^1(\aT)$. The function $\al_u$ has the following properties:
 \begin{enumerate}[(i)]
\item $\|\al_u\|_{L^1} = 1$.
\item $\al_u(e^{i\th}) = \il {\aD} P(z, e^{i\th})\,\Dl u(z)$.
\item $\al_u(e^{i\th})$ is lower semicontinuous.
\item $\al_u(e^{i\th}) \ge c > 0$ on $\aT$.
\item $\al_u(e^{i\th})$ need not to be necessarily bounded.
 \end{enumerate}

\section{Norms on Hardy Spaces}\label{norm}
\par Let us denote by $h^p_u(\aD),\, p > 1$, $u\in\E$, the space of harmonic functions $h $ on $\aD$ such that
\[ \|h\|^p_{u,p} = \varlimsup_{r\to 0^-} \int_{S_{u,r}} |h|^p\, d\mu_{u,r} < \infty.\]
By Corollary 3.2 in \cite{PS}, $h_u^p(\aD)\sbs h^p(\aD)$. Thus if $h\in h^p_u(\aD)$, then $h$ has radial boundary values $h^*$ on $\aT$. From the classical theory it is known that $h^*\in L^p(\aT)$ and $\|h\|_{h^p} = \|h^*\|_{L^p}$. We will see in the following theorem that the same holds for the functions in the new space.
\bT \label{hbvalue}Let $h\in h^p_u(\aD),\, p > 1$. Then $h^*\in L^p_u(\aT) := L^p(\aT, \mu_u)$ and
$ \|h\|_{u,p} = \|h^*\|_{L^p_u}$. \eT

\begin{proof} The least harmonic majorant on $\aD$ of the subharmonic function $|h|^p$ is the Poisson integral of $|h^*|^p$. By the Riesz Decomposition Theorem
\[|h(w)|^p=\il {\aT}|h^*(e^{i\th})|^pP(w,e^{i\th})\,d\lm(\th) +
\il {\aD}G(w,z)\Dl|h|^p(z),\]  where $P$ is the Poisson kernel and $G$ is the Green kernel.
\par By Lelong--Jensen formula and the monotone convergence theorem we have
\[\|h\|^p_{u,p}=\il {\aD}|h|^p\Dl u-\il {\aD}u\Dl|h|^p.\]
Again by the Riesz formula,
\begin{equation}\label{rf} u(z)=\il {\aD}G(z,w)\Dl u(w). \end{equation}
Hence, by Fubini--Tonnelli's  Theorem and the symmetry of the Green kernel
\[\il {\aD}u(z)\Dl|h|^p(z)=
\il {\aD}\left(\il {\aD}G(w,z)\Dl|h|^p(z)\right)\Dl u(w)\]
and
\begin{equation}\begin{aligned}
\|h\|^p_{u,p} &=\il {\aD}\left(|h(w)|^p-\il {\aD}G(w,z)\Dl|h|^p(z)\right)\Dl u(w)\notag\\
&=\il {\aD}\left(\il{\aT}|h^*(e^{i\th})|^pP(w,e^{i\th})\,d\lm(\th)\right)\Dl u(w)\\
&=\il{\aT}|h^*(e^{i\th})|^p\left(\il{\aD}P(w,e^{i\th})\Dl u(w)\right)\,d\lm (\th)\notag\\
&=\il{\aT}|h^*(e^{i\th})|^p\al_u(e^{i\th})\,d\lm (\th)\notag\\
&=\|h^*\|^p_{L^p_u}.
\end{aligned}\end{equation}
\end{proof}

This theorem has been proved also for the functions in $H^p_u(\aD)$ in \cite{MN} when $p > 1$, in \cite{Sibel} when $p \ge1$ and in \cite{Sh} when $p > 0$. We mention it here for the sake of completeness.
\bT\label{bvalue}
Let $f\in H^p(\aD),\,p > 0 $. Then $f\in H^p_u(\aD)$ if and only if $f^*\in L^p_u(\aT)$. Moreover, $\|f\|_{H^p_u} = \|f^*\|_{L^p_u}.$
\eT

\par In the proof of the theorem above we have deduced the norm of the functions $h\in h^p_u(\aD),\, p>1,$ to
\[\|h\|^p_{u,p} = \il {\aT}\left(\il \aD P(w, e^{i\th})\,\Dl u(w)\right)|h^*(e^{i\th})|^p\,d\lm.\]
Since the directional derivative of the Green kernel in the direction of the outward unit normal vector to $\aT$ is the Poisson kernel, i.e. $\frac{\bd }{\bd n}G(z, w)|_{z=e^{i\th}}=P(w, e^{i\th})$ ,  from the Riesz formula (\ref{rf})  we get
\[ \frac{\bd u}{\bd n}(e^{i\th}) = \il \aD P(w, e^{i\th})\,\Dl u(w)\]
and therefore the norm can be written as
\[\|h\|^p_{u,p} = \il {\aT}\frac{\bd u}{\bd n}(e^{i\th})|h^*(e^{i\th})|^p\,d\lm.\]
From this deduction it is clear that if $u\in\E$ is such that $\frac{\bd u}{\bd n}(e^{i\th})$ is bounded then $h^p_u(\aD) = h^p(\aD),\,p>1$.

\par For $u\in \E$, define $\E_u = \{ v\in \E: bv \le u \le b^{-1}v \text{ for some constant $b > 0$ near $\aT$}\}.$ It has been discussed in \cite{PS} that all the exhaustions in $\E_u$ generate the same weighted Hardy space $H^p_u(\aD)$ with the equivalent norms. We want to look into whether there is a canonical exhaustion in $\E_u$ determining the canonical norm on $H^p_u(\aD)$. Let $\E_0$ denote the set of all $u\in \E$ such that $\Dl u$ has a compact support in $\aD$. Then all the exhaustions in $\E_0$ generate the same space, the classical Hardy space $H^p(\aD)$, with the distinct norms and this is the largest space in our class.
\par For $u\in \E$ define
\begin{align*}
B_{u,p}(R)  & =\{f\in H^p_u(\aD): \|f\|_{H^p_u}\le R\} \text{ and }\\
B_\infty(R) & = \{f\in H^\infty(\aD):|f| \le R\}.
\end{align*}

\bT For $p > 0$, \[\bigcap_{u\in \E_0}B_{u,p}(1) = B_\infty(1).\]\eT 
\begin{proof}
The inclusion $B_\infty(1) \sbs \bigcap_{u\in \E_0}B_{u,p}(1)$ is clear. For the other way around, let $f\in H^\infty(\aD)\setminus B_\infty(1)$. Since $|f^*|^p \in L^1(\aT)$, by the Fatou's theorem
\[\int_{\aT}P(re^{i\hvi},e^{i\th})|f^*(e^{i\th})|^p\,d\lm \to |f^*(e^{i\hvi})|^p\]
$\lm$-a.e. on $\aT$. Hence there exists  $A\sbs\aT$  with  $\lm(A) >0$ such that
\begin{itemize}
\item $|f^*(e^{i\hvi})|> 1$  and
\item $\int_{\aT}P(re^{i\hvi},e^{i\th})|f^*(e^{i\th})|^p\,d\lm \to |f^*(e^{i\hvi})|^p$
\end{itemize}
for every $e^{i\hvi} \in A$. We may suppose that $1\in A$.

\par Since $u(z) = \int_{\aD} G(z, w)\,\Dl u(w)$, where $G(z,w)$ is the Green's function for the unit disk, and $\frac{\bd }{\bd n}G(z,w)|_{z=e^{i\th}}=P(w, e^{i\th})$,
\[\frac{\bd u}{\bd n}(e^{i\th}) = \int_\aD P(w, e^{i\th})\,\Dl u(w) = \al_u(e^{i\th}).\]
Also we have for $f\in H^p_u(\aD)$,
\[\|f\|^p_{H^p_u} = \int_{\aT}\frac{\bd u}{\bd n}(e^{i\th}) |f^*(e^{i\th})|^p\,d\lm.\]

\par Let $t_k\nea 1$ and $u_k(z) = G(z, t_k)$.  Then
\begin{equation*}\begin{aligned}
\|f\|^p_{H^p_{u_k}} & = \int_{\aT}P(t_k, e^{i\th})|f^*(e^{i\th})|^p\,d\lm\\
&\longrightarrow |f^*(1)|^p
\end{aligned}\end{equation*}
as $ k\to \infty $ because $1\in A$. Hence $f\not\in \bigcap_{u\in\E_0} B_{u,p}(1).$ The theorem follows.
\end{proof}
\par The above theorem suggests that it is impossible to select a canonical exhaustion determining the canonical norm on $H^p_u(\aD)$. 
\section{Weak-$*$ convergence of measures $\mu_{u,r}$}\label{wc}
\par While functions in $h^p_u(\aD)$, $p>1$, have radial limits $\mu_u$-a.e., we are interested in the analogs of more subtle classical properties of boundary values. For example, if $h\in h^p(\aD)$ then it is known that the measures $h(re^{i\th})\lm(\th)$ converge weak-$*$ in $C^*(\aT)$ to $h^*(e^{i\th})\lm(\th)$ as $r\to1^-$.
\par In this section we will establish the analog of this statement.
\bT \label{wlimit} Let $h \in h^p_u(\aD),\, p > 1$. Then the measures $\{h\mu_{u,r}\}$ converge weak-${\ast}$ to $h^*\mu_u$ in $C^*(\oD)$ when $r\to0^-$. \eT
\begin{proof}
Since the space $C(\oD)$ is separable, the weak-$*$ topology on the balls in $C^*(\oD)$ is metrizable. Thus it suffices to show that for any sequence $r_j\nea0$ and any $\phi\in C(\oD)$ we have
\[\lim_{j\to\infty}\il{S_{u,r_j}}\phi h\,d\mu_{u,r_j}=\il {\aT}\phi h^*\,d\mu_u.\]
\par We introduce functions
\[p_r(e^{i\th})=\il{S_{u,r}}P(z,e^{i\th})\,d\mu_{u,r}(z)=\il{B_{u,r}}P(z,e^{i\th})\,\Dl u(z),\] where the second equality is due to the Lelong--Jensen formula (Theorem \ref{jlf}).  Hence $p_r(e^{i\th})\nea\al_u(e^{i\th})$.
\par  Due to the uniform continuity of $\phi$ and the formula for $P(z,e^{i\th})$,  for every $\th\in[0,2\pi]$ and for every $\eps>0$ there is $\dl>0$ such that $|P(z,e^{i\th})|<\eps$ when $z$ is close to boundary and $|z-e^{i\th}|>\dl$ and $|\phi(z)-\phi(e^{i\th})|<\eps$ when $|z-e^{i\th}|\le\dl$. Hence, when $r$ is sufficiently close to 0,
\begin{equation}\begin{aligned}
&\left|\il { S_{u,r}}\phi(z)P(z,e^{i\th})\,d\mu_{u,r}(z)-
\il{S_{u,r}}\phi(e^{i\th})P(z,e^{i\th})\,d\mu_{u,r}(z)\right| \notag\\
\le&\il{S_{u,r}\sm\oD(e^{i\th},\dl)}|\phi(z)-\phi(e^{i\th})|P(z,e^{i\th})\,d\mu_{u,r}(z)\notag\\
& +\il{S_{u,r}\cap\oD(e^{i\th},\dl)}|\phi(z)-\phi(e^{i\th})|P(z,e^{i\th})\,d\mu_{u,r}(z) \notag\\
\le&2M\eps+\eps p_r(e^{i\th}),\notag
\end{aligned}\end{equation}
where $\aD(e^{i\th}, \delta)$ is the disk of radius $\delta$ and center at $e^{i\th}$ and $M$ is the uniform norm of $\phi$ on $\ovr{\aD}$.
\par Now,
\begin{equation}\begin{aligned}\il{S_{u,r}}\phi(z)h(z)\,d\mu_{u,r}(z)=&
\il{S_{u,r}}\phi(z)\left(\il{\aT}h^*(e^{i\th})P(z,e^{i\th})\,d\lm(\th)\right)\,d\mu_{u,r}(z)\notag\\=
&\il{\aT}h^*(e^{i\th})\left(\il{S_{u,r}}\phi(z)P(z,e^{i\th})\,d\mu_{u,r}(z)\right)\,d\lm(\th).\notag
\end{aligned}\end{equation}
Hence,
\begin{equation}\begin{aligned}
&\left|\il{S_{u,r}}\phi(z)h(z)\,d\mu_{u,r}(z)-
\il{\aT}\phi(e^{i\th})h^*(e^{i\th})\,d\mu_u(\th)\right|\notag\\
\le& \left|\il{S_{u,r}}\phi(z)h(z)\,d\mu_{u,r}(z)- \il {\aT} \phi(e^{i\th})h^*(e^{i\th})p_r(e^{i\th})\,d\lm(\th)\right|\notag\\
&+\left|\il {\aT} \phi(e^{i\th})h^*(e^{i\th})p_r(e^{i\th})\,d\lm(\th)-
\il{\aT}\phi(e^{i\th})h^*(e^{i\th})\,d\mu_u(\th)\right|\notag\\
=&\left|\il{\aT} h^*(e^{i\th})\left(\il{S_{u,r}}(\phi(z)-
\phi(e^{i\th}))P(z,e^{i\th})\,d\mu_{u,r}(z)\right)\,d\lm(\th)\right|\notag\\
&+\left|\il {\aT} \phi(e^{i\th})h^*(e^{i\th})\left(p_r(e^{i\th})-\al_u(e^{i\th})\right)\,d\lm(\th)\right|\notag\\
\le&\eps\il{\aT}\left |h^*(e^{i\th})\right|(2M+p_r(e^{i\th}))\,d\lm(\th)  + M\il {\aT} \left|h^*(e^{i\th})\right|\left|p_r(e^{i\th})-\al_u(e^{i\th})\right|\,d\lm(\th).\notag
\end{aligned}\end{equation}
 Now,
\begin{equation}\begin{aligned}
 \il{\aT}\left |h^*(e^{i\th})\right|(2M+p_r(e^{i\th}))\,d\lm(\th) \le &\il{\aT}\left |h^*(e^{i\th})\right|(2M+\al_u(e^{i\th}))\,d\lm(\th)\notag\\
\le & 2M\|h^*\|_{L^p} + \|h\|_{u,p}.\notag
\end{aligned}\end{equation}
Since $\left|p_r(e^{i\th})-\al_u(e^{i\th})\right| \sea 0$ and $\left|p_r(e^{i\th})-\al_u(e^{i\th})\right| < \al_u(e^{i\th})$ with $\left|h^*(e^{i\th})\right|\al_u(e^{i\th}) \in L^1(\aT)$, by the monotone convergence theorem,
\[\il {\aT} \left|h^*(e^{i\th})\right|\left|p_r(e^{i\th})-\al_u(e^{i\th})\right|\,d\lm(\th) \to 0\]
Thus, since $\eps$ is arbitraty,
\[\left|\il{S_{u,r}}\phi(z)h(z)\,d\mu_{u,r}(z)-
\il{\aT}\phi(e^{i\th})h^*(e^{i\th})\,d\mu_u(\th)\right|\to 0.\]
The proof is complete.
\end{proof}
\bC If $p > 1$, the measures $\mu_{u,r}$ converge weak-$*$ to $\mu_u$ in the dual of $h^p_u(\aD)$ when $r \to 0^-$.  \eC
\begin{proof} 
For $\phi\in C(\ovr{\aD})$, from the theorem above we have
\[\lim_{r\to 0^-}\int_{S_{u,r}}\phi h\,d\mu_{u,r} = \int_{\aT}\phi h^*\,d\mu_u\]
for every $h\in h^p_u(\aD)$. In particular, if we take $\phi \equiv 1$ we get
\[\lim_{r\to 0^-}\int_{S_{u,r}} h\,d\mu_{u,r} = \int_{\aT} h^*\,d\mu_u\]
for every $h\in h^p_u(\aD)$. The corollary follows. 
\end{proof}
\par The corollary above improves the result of Demailly that the measures $\mu_{u,r}$ converge weak-$*$ to $\mu_u$ in $C^*(\ovr{\aD})$. The convergence is much stronger indeed. 
\section{Convergence of dilations}\label{cd}
\par For $0 < t < 1$ we define the dilations of a function $f$ as $f_t(z) = f(tz)$. In the classical theory, if $f\in H^p(\aD)$ then the $H^p$-norm of $f$ is defined to be the limit of the $L^p$-norms of the dilations as $t\to 1$, that is, 
\begin{equation} \label{cnorm} \|f\|_{H^p} = \lim_{t\to 1}\|f_t\|_{L^p} = \lim_{t\to 1}\left(\int_\aT |f(te^{i\th}|^p\,d\lm\right)^{1/p}. \end{equation}
Moreover, the limit on the right converges to the $L^p$-norm of the boundary value function $f^*$ as one could expect. Thus we have 
\begin{equation} \label{enorm} \|f\|_{H^p} = \|f^*\|_{L^p}.\end{equation}
Also the dilations converge to $f^*$ in the $L^p$-norm and so do converge to $f$ in the $H^p$-norm, that is, 
\begin{equation}\label{nconv}\lim_{t\to 1}\|f_t - f^*\|_{L^p} = \lim_{t\to1}\|f_t - f\|_{H^p} = 0.\end{equation}
Note that here the term ``norm'' has been abused to indicate $\|\cdot\|_H^p$ or $\|\cdot\|_{L^p}$ for all $p > 0$. 

\par We want to establish the analog of these statements to the new theory.  The analog of (\ref{enorm}) has been discussed in section \ref{norm}. We will establish the analog of (\ref{cnorm}) and (\ref{nconv}) in this section. 

\bL\label{Poisson} For $0 < t < 1,$ \[\int_{\aT}P(te^{i\hvi}, e^{i\th})P(z, e^{i\th})\,d\th = P(tz, e^{i\hvi}).\]\eL
 \begin{proof}
If $z=re^{i\psi}$ then the Poisson kernel  is
\[P(z,e^{i\th})=P_r(\psi-\th)=\sum_{k=-\infty}^{\infty}r^{|k|}e^{ik(\psi-\th)}.\]
If $0\le r<1$ this series converges uniformly in $\th$ and
\[\il{\aT}P(te^{i\phi},e^{i\th})P(z,e^{i\th})\,d\lm(\th)=
\sum_{k=-\infty}^{\infty}r^{|k|}e^{ik\psi}\il{\aT}P(te^{i\phi},e^{i\th})e^{-ik\th}\,d\lm(\th).\]
If $k<0$ then the function $e^{-ik\th}$ is the boundary value of $\zeta^{|k|}$ where $\zeta=te^{i\th}$ while when $k\ge0$ its is the boundary value of $\bar{\zeta}^{|k|}$ where $\zeta=te^{i\th}$. Hence
\[\il{\aT}P(te^{i\phi},e^{i\th})e^{-ik\th}\,d\lm(\th)=t^{|k|}e^{-ik\phi}\] when $k<0$ and
\[\il{\aT}P(te^{i\phi},e^{i\th})e^{-ik\th}\,d\lm(\th)=t^{|k|}e^{-ik\phi}\] when $k\ge0$. Thus
\[\sum_{k=-\infty}^{\infty}r^{|k|}e^{ik\psi}\il{\aT}P(te^{i\phi},e^{i\th})e^{-ik\th}\,d\lm(\th)=
\sum_{k=-\infty}^{\infty}r^{|k|}t^{|k|}e^{ik(\psi-\phi)}\] and we see that
\[\il{\aT}P(te^{i\phi},e^{i\th})P(z,e^{i\th})\,d\lm(\th)=P(tz,e^{i\phi}).\]
\end{proof}

\bT\label{approx}
Let $f\in H^p_u(\aD),\, p > 0 $. Then we have
\begin{enumerate}[(i)]
\item $ \lim_{t\to 1} \|f_t\|_{H^p_u} =\|f\|_{H^p_u}$\,\,\, and
\item $ \lim_{t\to 1} \|f_t-f\|_{H^p_u} =0$.
\end{enumerate}
\eT

\begin{proof}
By Theorem \ref{bvalue} and Fatou's lemma,
\begin{equation}\label{E1} \|f\|_{H^p_u}^p = \|f^*\|^p_{L^p_u}\le \liminf_{t\to1}\int_{\aT}|f(te^{i\th})|^p\,d\mu_u.\end{equation}
Recall that \[\al_u(e^{i\th}) = \lim_{r\to0}\int_{S_{u,r}}P(z, e^{i\th})\,d\mu_{u,r}(z)\] and the integral on the right hand side is an increasing function of $r$.
Therefore,
\begin{align*}
\int_{\aT}|f(te^{i\th})|^p\,d\mu_u(\th) & = \int_{\aT}|f(te^{i\th})|^p\left(\lim_{r\to0}\int_{S_{u,r}}P(z, e^{i\th})\,d\mu_{u,r}(z)\right)\,d\lm(\th)\\
&= \lim_{r\to 0}\int_{\aT}|f(te^{i\th})|^p\left(\int_{S_{u,r}}P(z, e^{i\th})\,d\mu_{u,r}(z)\right)\,d\lm(\th).
\end{align*}
Given $\eps>0$, there exists $r_0<0$ such that for $r_0<r<0$,
\begin{align*}
\int_{\aT}|f(te^{i\th})|^p\,d\mu_u(\th) -\eps &  \le \int_{\aT}|f(te^{i\th})|^p\left(\int_{S_{u,r}}P(z, e^{i\th})\,d\mu_{u,r}(z)\right)\,d\lm(\th)\\
 & = \int_{S_{u,r}}\left(\int_{\aT}|f(te^{i\th})|^pP(z, e^{i\th})\,d\lm(\th)\right)\,d\mu_{u,r}(z).
\end{align*}
Using the subharmonicity of $|f|^p$ we can write
\[ |f(te^{i\th})|^p \le \int_{\aT}|f^*(e^{i\hvi})|^pP(te^{i\th}, e^{i\hvi})\,d\lm(\hvi)\]
and
\begin{align*}
\int_{\aT}|f(te^{i\th})|^pP(z, e^{i\th})\,d\lm(\th) & \le \int_{\aT}\left(\int_{\aT}|f^*(e^{i\hvi})|^pP(te^{i\th}, e^{i\hvi})\,d\lm(\hvi)\right) P(z, e^{i\th})\,d\lm(\th)\\
& = \int_{\aT}|f^*(e^{i\hvi})|^p\left(\int_{\aT}P(te^{i\hvi},e^{i\th})P(z,e^{i\th})\,d\lm(\th)\right)\,d\lm(\hvi)\\
&=\int_{\aT}|f^*(e^{i\hvi})|^pP(tz,e^{i\hvi})\,d\lm(\hvi).
\end{align*}
For fixed $0 > r > r_0$,
\begin{align*}
\lim_{t\to 1} \int_{\aT}|f(te^{i\th})|^p\,d\mu_u(\th) -\eps &\le \lim_{t\to 1} \int_{S_{u,r}}\left(\int_{\aT}|f^*(e^{i\hvi})|^pP(tz,e^{i\hvi})\,d\lm(\hvi)\right)\,d\mu_{u,r}(z)\\
& =\lim_{t\to 1}\int_{\aT}|f^*(e^{i\hvi})|^p\left(\int_{S_{u,r}}P(tz,e^{i\hvi})\,d\mu_{u,r}(z)\right)\,d\lm(\hvi)\\
& =\int_{\aT}|f^*(e^{i\hvi})|^p\left(\int_{S_{u,r}}P(z,e^{i\hvi})\,d\mu_{u,r}(z)\right)\,d\lm(\hvi).
\end{align*}
 Hence
\begin{equation}\begin{aligned}\label{E2}
& \lim_{t\to 1} \int_{\aT}|f(te^{i\th})|^p\,d\mu_u(\th) -\eps \\
& \le \lim_{r\to 0}\int_{\aT}|f^*(e^{i\hvi})|^p\left(\int_{S_{u,r}}P(z, e^{i\hvi})\,d\mu_{u,r}(z)\right)\,d\lm(\hvi)\\
& = \int_{\aT}|f^*(e^{i\hvi})|^p \al_u(e^{i\hvi})\,d\lm(\hvi)\\
& =  \|f\|^p_{H^p_u}.
\end{aligned}\end{equation}
Thus from (\ref{E1}) and (\ref{E2}) we have $\lim_{t\to 1}\|f_t\|_{H^p_u}=\|f\|_{H^p_u}$.

\par To prove $(ii)$, let $\eps > 0$ be given. Since $f^*\in L^p_u(\aT)$, there exists $\dl> 0$ such that for any set $E\sbs \aT$ with $\mu_u(E)< \dl$ we have
\[\int_E |f^*|^p\,d\mu_u < \frac{\eps}{2}.\]
Also, $f_t\to f^*$\, $\mu_u$-a.e. on $\aT$. Apply Egorov's theorem to get a set $E\sbs \aT$ such that $\mu_u(E) < \dl$ and $f_t \to f^*$ uniformly on $\aT\sm E$. Hence
\[\lim_{t\to 1}\int_{\aT\sm E}|f_t|^p\,d\mu_u = \int_{\aT\sm E}|f^*|^p\,d\mu_u.\]
From part (1) we have
\[\lim_{t\to 1}\int_{\aT}|f_t|^p\,d\mu_u = \int_{\aT}|f^*|^p\,d\mu_u.\]
Combining them we get
\[\lim_{t\to 1}\int_E |f_t|^p\,d\mu_u = \int_E |f^*|^p\,d\mu_u,\]
the right hand side of which is less than $\eps/2$. Therefore
\[\int_E |f_t|^p\,d\mu_u < \eps\]
for $t$ near $1$. Now
\[\int_{\aT} |f_t - f^*|^p\,d\mu_u = \int_{\aT\sm E} |f_t - f^*|^p\,d\mu_u + \int_E |f_t - f^*|^p\,d\mu_u.\]
Since $f_t \to f^*$ uniformly on $\aT\sm E$,
\[\int_{\aT\sm E} |f_t - f^*|^p\,d\mu_u  \to 0\]
and for $t$ near $1$,
\begin{eqnarray*}
\int_E |f_t - f^*|^p\,d\mu_u  & \le  & 2^p\left(\int_E |f_t|^p\,d\mu_u + \int_E |f^*|^p\,d\mu_u\right)\\
& <  & 2^{p+1}\eps.
\end{eqnarray*}
Therefore
\[\int_{\aT} |f_t - f^*|^p\,d\mu_u \to 0\]
and $(2)$ follows once again by Theorem \ref{bvalue}.
\end{proof}

\par As a consequence of this theorem we have the following result about density.

\bC\label{density}
Polynomials are dense in $H^p_u(\aD),\, p > 0$.
\eC

\end{document}